\documentclass[aap,MSNbibl,citesort,dvips]{arximspdf}

\doi{10.1214/11-AAP837} 
\volume{22}
\issue{6}
\pubyear{2012}
\firstpage{2429}
\lastpage{2459}

\makeatletter

\renewcommand{\widehat}{\hat}

\newcommand{\acal}{\mathcal{A}}
\newcommand{\bcal}{\mathcal{B}}
\newcommand{\dcal}{\mathcal{D}}
\newcommand{\ecal}{\mathcal{E}}
\newcommand{\fcal}{\mathcal{F}}
\newcommand{\ical}{\mathcal{I}}
\newcommand{\kcal}{\mathcal{K}}
\newcommand{\pcal}{\mathcal{P}}
\newcommand{\rcal}{\mathcal{R}}
\newcommand{\scal}{\mathcal{S}}
\newcommand{\tcal}{\mathcal{T}}
\newcommand{\ucal}{\mathcal{U}}

\newcommand{\real}{\mathbb{R}}
\newcommand{\A}{\acal}
\renewcommand{\l}{\lambda}

\newcommand{\eps}{\varepsilon}

\newcommand{\bfz}{\mathbf{z}}

\newcommand{\prob}{\mathbb{P}}
\newcommand{\expec}{\mathbb{E}}
\newcommand{\E}{\mathbb{E}}
\renewcommand{\P}{\mathbb{P}}

\newcommand{\poly}{\operatorname{poly}}

\newcommand{\weight}{w}
\newcommand{\dist}{d}
\newcommand{\nsqneq}{[n]^2_{\neq}}
\newcommand{\nfourneq}{[n]^4_{\neq}}
\renewcommand{\path}{\mathrm{Path}}
\newcommand{\short}{\mbox{\L}^-_2}
\newcommand{\notshort}{\mbox{\L}^+_2}

\newcommand{\mixed}{\mathbb{T}}
\newcommand{\phy}{\mathbb{Y}}
\newcommand{\minnu}{\underline{\nu}}
\newcommand{\gmassump}{\Theta\mbox{-}\mathrm{M}}
\newcommand{\tv}[1]{\|#1\|_{\mathrm{TV}}}
\newcommand{\qmin}{Q_{\mathrm{min}}}
\newcommand{\hatupsilon}{\widehat{\Upsilon}}
\newcommand{\hatucal}{\widehat{\ucal}}
\newcommand{\hatkcal}{\widehat{\kcal}}

\newtheorem{theorem}{Theorem}

\newproclaim{definition}{Definition}[section]

\newtheorem{lemma}{Lemma}[section]

\newproclaim{remark}{Remark}[section]

\newproclaim{assumption}{Assumption}

\makeatother

\begin{document}
\begin{frontmatter}

\title{Phylogenetic mixtures:
Concentration of measure in the large-tree limit}
\runtitle{Phylogenetic mixtures in the large-tree limit}

\begin{aug}
\author[A]{\fnms{Elchanan} \snm{Mossel}\thanksref{t1}\ead[label=e1]{mossel@stat.berkeley.edu}}
\and
\author[B]{\fnms{Sebastien} \snm{Roch}\corref{}\thanksref{t2}\ead[label=e2]{roch@math.ucla.edu}}
\runauthor{E. Mossel and S. Roch}
\affiliation{University of California, Berkeley, and University~of~California,~Los~Angeles}
\address[A]{Departments of Statistics\\
\quad and Computer Science\\
University of California\\
Berkeley, California 94720\\
USA\\
\printead{e1}}
\address[B]{Department of Mathematics\\
\quad and Bioinformatics Program\\
University of California\\
Los Angeles, California 90095\\
USA\\
\printead{e2}} 
\end{aug}

\thankstext{t1}{Supported by DMS 0548249 (CAREER) award,
by DOD ONR Grant N000141110140, by ISF
Grant 1300/08 and by ERC Grant PIRG04-GA-2008-239137.}

\thankstext{t2}{Supported by NSF Grant DMS-10-07144.}

\received{\smonth{8} \syear{2011}}
\revised{\smonth{12} \syear{2011}}

%
\begin{abstract}
The reconstruction of phylogenies from DNA or protein sequences is a
major task of computational evolutionary biology. Common phenomena,
notably variations in mutation rates across genomes and incongruences
between gene lineage histories, often make it necessary to model
molecular data as originating from a mixture of phylogenies. Such mixed
models play an increasingly important role in practice.

Using concentration of measure techniques, we show that mixtures of
large trees are typically identifiable. We also derive sequence-length
requirements for high-probability reconstruction.
\end{abstract}

%

\end{frontmatter}

\section{Introduction}\label{sectionintroduction}

Phylogenetics~\cite{SempleSteel03,Felsenstein04} is centered around
the reconstruction of evolutionary histories
from molecular data extracted from modern species.
The assumption is that molecular data consists of aligned sequences and
that each position in the sequences evolves independently according to a
Markov model on a tree, where the key parameters are
(see Section~\ref{sectionpreliminaries} for formal definitions):
\begin{itemize}
\item \textit{Rate matrix.}
An $r \times r$ mutation rate matrix $Q$,
where $r$ is the alphabet size.
A~typical alphabet is the set of nucleotides $\{\mathrm{A},\mathrm
{C},\mathrm{G},\mathrm{T}\}$,
but here we allow more general state spaces.
Without loss of generality, we denote the alphabet by
$\rcal= [r] = \{1,\ldots,r\}$.
The $(i,j)$th entry of $Q$ encodes the rate at which state $i$ mutates
into state $j$.
\item \textit{Binary tree.}
An evolutionary tree $T$, where the leaves are the modern species and
each branching represents a past speciation event.
The leaves are labeled with names of species.
Without loss of generality, we assume
the labels are $X = [n]$.
\item \textit{Branch lengths.}
For each edge $e$,
we have a scalar branch length $\weight_e$
which measures the expected total number of substitutions per site
along edge~$e$.
Roughly speaking, $\weight_e$ is the amount of
mutational change between the end points of $e$.
\end{itemize}
The classical problem in phylogenetics can be
stated as follows:
\begin{itemize}
\item
\textit{Phylogenetic tree reconstruction} (\textit{PTR}):
\textit{Unmixed case.}
Given $n$ molecular sequences of length $k$,
\[
\{s_a = (s^i_a)_{i=1}^k\}_{a\in[n]}
\]
with $s^i_a \in[r]$,
which have evolved according to the process above with independent sites,
reconstruct the topology of the evolutionary tree.
\end{itemize}
There exists a vast theoretical literature on this
problem; see, for example,~\cite{SempleSteel03} and
references therein.

However,
various
phenomena, notably variations in mutation rates across
genomes
and incongruences between gene lineage histories, often
make it necessary to model molecular data as originating from a
\textit{mixture}
of different phylogenies.

Here, using concentration of measure techniques,
we show that mixtures of large trees are
typically identifiable. By typically, we mean
informally that our results hold
under conditions guaranteeing that
the tree topologies present in the mixture are sufficiently
distinct. (See Section~\ref{secresults} for
a careful statement of the theorems.)
In particular, we give a broad new class of conditions
under which mixtures are identifiable, and
we extend, to more general substitution
models, previous results on the total variation
distance between Markov models on trees.
Our proofs are constructive in that we provide
a computationally efficient reconstruction algorithm.
We also derive sequence-length requirements for
high-probability reconstruction.

Our
identifiability and reconstruction
results represent an important first step toward dealing
with more biologically relevant mixture models
(such as the ones mentioned above) in which
the tree topologies tend to be similar.
In particular, in a recent related
paper~\cite{MosselRoch11a}, we have used
the techniques developed here to
reconstruct common rates-across-sites models.

\subsection{Related work}\label{sectionrelated}

Most prior theoretical work on mixture models
has focused on the question of \textit{identifiability}.
A class of phylogenetic models is identifiable
if any two models in the class produce different
data distributions.
It is well known that unmixed phylogenetic models
are typically identifiable~\cite{Chang96}.
\textit{This is not the case in general for mixtures of phylogenies.}
For instance, Steel et al.~\cite{StSzHe94} showed that for any two trees
one can find a random scaling on each of them, such that their data
distributions
are identical.\vadjust{\goodbreak} Hence it is hopeless, in general, to reconstruct
phylogenies under mixture models.
See also \cite
{EvansWarnow04,MatsenSteel07,MaMoSt08,StefankovicVigoda07a,StefankovicVigoda07b,Steel09}
for further examples of this type.\looseness=-1

However, the negative examples constructed in the references above are
not necessarily typical.
They use special features of the mutation models (and their invariants)
and allow themselves quite a bit of flexibility in setting up the
topologies and branch lengths.
In fact, recently a variety of more standard mixture models have been
shown to be
identifiable. These include the common
$\mbox{GTR}+\Gamma$ model~\cite{AlAnRh08,WuSusko10} and
$\mbox{GTR}+\Gamma+\mbox{I}$ model~\cite{ChaiHousworth11},
as well as some covarion models~\cite{AllmanRhodes06},
some group-based models~\cite{AlPeRhSu11}
and so-called $r$-component identical tree mixtures~\cite{RhodesSullivant10}.
Although these results do not provide practical algorithms for reconstructing
the corresponding mixtures, they do give hope that these problems
may be tackled successfully.

Beyond the identifiability question, there seems to have been little rigorous
work on reconstructing phylogenetic mixture models. One positive
result is the case of the molecular clock assumption with across-sites
rate variation~\cite{StSzHe94}, although no sequence-length
requirements are provided.
There is a large body of work on practical reconstruction
algorithms for various types of mixtures, notably rates-across-sites models
and covarion-type models, using mostly likelihood and Bayesian methods;
see, for example,~\cite{Felsenstein04} for references.
But the optimization problems they attempt to solve are
likely NP-hard~\cite{ChorTuller06,Roch06}.
There also exist many techniques for testing for the presence of a mixture
(e.g., for testing for rate heterogeneity), but such tests typically require
the knowledge of the phylogeny; see, for example,~\cite{HuelsenbeckRannala97}.

Here we give both identifiability and reconstruction
results. The proof of our main results relies
on the construction of a \textit{clustering statistic} that
discriminates between distinct phylogenies.
A similar approach was
used recently in~\cite{MosselRoch11a}.
There, however,
the problem was to distinguish between phylogenies
with the \textit{same} topology,
but \textit{different} branch lengths.
In the current work, a main technical challenge
is to analyze the simultaneous behavior
of such a clustering statistic on
\textit{distinct} topologies.
A similar statistic was also used
in~\cite{SteelSzekely06} to prove a special case
of Theorem~\ref{thm2} below. However,
in contrast to~\cite{SteelSzekely06}, our main result
requires that
a clustering statistic be constructed
based only on data generated by
the mixture---that is, \textit{without} prior
knowledge of the topologies to be
distinguished. Finally, unlike~\cite{MosselRoch11a}
and~\cite{SteelSzekely06},
we consider the more general
GTR model.

\section{Definitions and results}

\subsection{Basic definitions}

\subsubsection*{Phylogenies}
A phylogeny is a graphical representation
of the speciation history of a group of organisms.
The leaves typically correspond to current species.
Each branching indicates a speciation event.
Moreover we associate to each edge a\vadjust{\goodbreak}
positive weight.
This weight can be thought roughly as the
amount of evolutionary change on the edge.
More formally, we make the following definitions; see, for example,
\cite{SempleSteel03}.
Fix a set of leaf labels $X = [n] = \{1,\ldots,n\}$.
%
\begin{definition}[(Phylogeny)]\label{defphylo}
A \textit{weighted binary phylogenetic $X$-tree}
(or \textit{phylogeny})
$T = (V,E;\phi;\weight)$ is a tree
with vertex set $V$,
edge set $E$,
leaf set $L$ with $|L| = n$,
and a bijective mapping $\phi\dvtx X \to L$
such that:
\begin{enumerate}
\item[(1)] The degree of all internal vertices
$V-L$ is exactly $3$.
\item[(2)] The edges are assigned
weights $\weight\dvtx E \to(0,+\infty)$.
\end{enumerate}
We let $\tcal_l[T] = (V,E;\phi)$ be the
\textit{leaf-labelled topology} of $T$.
\end{definition}
%
\begin{definition}[(Tree metric)]
A phylogeny $T = (V,E;\phi;\weight)$
is naturally equipped with a
\textit{tree metric}
$\dist_T\dvtx X\times X \to(0,+\infty)$ defined as follows:
\[
\forall a,b \in X\qquad \dist_T(a,b)
= \sum_{e\in\path_T(\phi(a),\phi(b))} \weight_e,
\]
where $\path_T(u,v)$ is the set of edges on the path between
$u$ and $v$ in $T$.
We will refer to $\dist_T(a,b)$ as the \textit{evolutionary
distance} between $a$ and $b$.
In a slight abuse of notation,
we also sometimes use
$\dist_T(u,v)$ to denote the
evolutionary distance as above between
any two vertices $u,v$ of $T$.
\end{definition}

We will restrict ourselves to the following standard special case.
%
\begin{definition}[(Regular phylogenies)]
Let $0 < f \leq g < +\infty$.
We denote by $\phy^{(n)}_{f,g}$
the set of phylogenies $T = (V,E; \phi;\weight)$
with $n$ leaves
such that $f \leq\weight_e \leq g$,
$\forall e\in E$. We also
let $\phy_{f,g} = \bigcup_{n \geq1}
\phy^{(n)}_{f,g}$.
\end{definition}

\subsubsection*{GTR model}
A commonly used model of DNA sequence evolution
is the following \textit{GTR model}; see, for example, \cite
{SempleSteel03}. We first define an appropriate
class of rate matrices.
%
\begin{definition}[(GTR rate matrix)]
Let $\rcal$ be a set of character states
with $r = |\rcal|$. Without loss of generality
we assume that $\rcal= [r]$.
Let $\pi$ be a probability distribution on
$\rcal$ satisfying $\pi_x > 0$ for all
$x \in\rcal$.
A \textit{general time-reversible} (\textit{GTR})
\textit{rate matrix} on $\rcal$, with
respect to stationary distribution $\pi$, is
an $r \times r$ real-valued matrix $Q$ such that:
\begin{enumerate}
\item[(1)] $Q_{xy} > 0$ for all $x\neq y \in\rcal$.

\item[(2)] $\sum_{y \in\rcal} Q_{x y} = 0$, for all $x \in\rcal$.

\item[(3)] $\pi_x Q_{x y} = \pi_y Q_{y x}$, for all $x, y \in\rcal$.
\end{enumerate}
By the reversibility assumption, $Q$
has $r$ real eigenvalues
\[
0 = \Lambda_1 > \Lambda_2 \geq\cdots\geq\Lambda_{r}.
\]
We normalize $Q$ by fixing $\Lambda_2 = -1$.\vadjust{\goodbreak}
\end{definition}
%
\begin{definition}[(GTR model)]\label{defmmt}
Consider the following stochastic
process.
We are given a phylogeny $T = (V,E;\phi;\weight)$
and a finite set $\rcal$ with $r$ elements.
Let $\pi$ be a probability distribution
on $\rcal$ and
$Q$ be a GTR rate matrix
with respect to $\pi$.
Associate to each edge $e\in E$ the stochastic matrix
\[
M(e) = \exp(\weight_e Q).
\]
The process runs as follows. Choose an arbitrary root $\rho\in V$.
Denote by $E_\downarrow$ the set $E$ directed away from the root. Pick
a state for the root at random according to~$\pi$. Moving away from the
root toward the leaves, apply the channel $M(e)$ to each edge $e$
independently. Denote the state so obtained $s_V = (s_v)_{v\in V}$. In
particular, $s_{L}$~is the state at the leaves, which we also denote by
$s_X$. More precisely, the joint distribution of $s_V$ is given by
\[
\mu_V(s_V) = \pi_{\rho}(s_\rho)
\prod_{e = (u,v) \in E_\downarrow}
[M(e)]_{s_{u} s_{v}}.
\]
For $W \subseteq V$, we denote by $\mu_W$ the marginal
of $\mu_V$ at $W$. Under this model, the weight $\weight_e$ is
the expected number of substitutions on edge $e$
in the
continuous-time process.
We denote by $\dcal[T,Q]$
the probability distribution of $s_V$.
We also let $\dcal_l[T,Q]$ denote
the probability distribution of
\[
s_X \equiv\bigl(s_{\phi(a)}\bigr)_{a\in X}.
\]
\end{definition}

More generally, we consider $k$ independent samples
$\{s^{i}_{V}\}_{i=1}^k$ from the model above, that is, $s^1_V,
\ldots,s^k_V$ are i.i.d. $\dcal[T,Q]$. We think of $(s_v^i)_{i=1}^k$ as
the sequence at node $v \in V$. Typically, $\rcal= \{\mathrm{A},
\mathrm{G},\mathrm{C},\mathrm{T}\}$ and the model describes how DNA
sequences stochastically evolve by point mutations along an
evolutionary tree under the assumption that each site in the sequences
evolves independently. When considering many samples
$\{s^{i}_{V}\}_{i=1}^k$, we drop the subscript to refer to a single
sample~$s_V$.

\subsubsection*{Mixed model}
We introduce the basic mixed model
which will be the focus of this paper.
We will use the following definition.
We assume that $Q$ is fixed
and known throughout.
%
\begin{remark}[(Unknown rate matrix)]
See the concluding remarks
for an extension of our techniques
when $Q$ is unknown.
\end{remark}
%
\begin{definition}[($\Theta$-mixture)]
Let $\Theta$ be a positive integer.
In the \textit{$\Theta$-mixture model},
we consider a finite set of
phylogenies
\[
\mixed=
\{T_\theta= (V_\theta,E_\theta;\phi_\theta;\weight_\theta)\}
_{\theta=
1}^{\Theta}
\]
on the same set of leaf labels $X = [n]$ and a positive probability
distribution $\nu= (\nu_{\theta})_{\theta=1}^\Theta$ on $[\Theta]$.
Consider $k$ i.i.d. random variables $N^1, \ldots, N^k$ with
distribution~$\nu$. Then, conditioned on $N^1, \ldots, N^k$, the
samples $\{s^{i}_X\}_{i=1}^k$ generated under the $\Theta$-mixture
model $(\mixed,\nu,Q)$ are independent with conditional distribution
$s^j_X \sim\dcal_l[T_{N^j}, Q]$, $j=1,\ldots,k$. We denote by
$\dcal_l[(\mixed,\nu,Q)]$ the probability distribution of $s^1_X$. We
will refer to $T_\theta$ as the \textit{$\theta$-component} of the
mixture $(\mixed,\nu,Q)$.
\end{definition}

We assume that $\Theta$ is fixed and known throughout.
As above, we drop the superscript to refer to a single
sample $s_X$ with corresponding component indicator
$N$. To simplify notation,
we let
\[
d_{T_{\theta}} = d_\theta\qquad \forall\theta\in[\Theta].
\]

\subsubsection*{Some notation}
We will\vspace*{1pt} use the notation $[n]^2 = \{(a,b)
\in[n]\times[n]\dvtx a \leq b\}$, $[n]^2_{=} = \{(a,a)\}_{a\in[n]}$ and
$\nsqneq= [n]^2 - [n]^2_{=}$. We also denote by $\nfourneq$ the set of
pairs $(a_1,b_1), (a_2,b_2)\in\nsqneq$ such that $(a_1,b_1)
\neq(a_2,b_2)$ (as pairs). We use the notation $\poly(n)$ to denote the
growth condition usually written $\Theta(n^C)$ for some $C > 0$.

\subsection{Main results}
\label{secresults}

We make the following assumptions on the mutation model.
%
\begin{assumption}\label{assumpgeneral}
Let
$0 < f \leq g < +\infty$,
and $\minnu> 0$.
We will use the following set of assumptions on
a $\Theta$-mixture model $(\mixed,\nu,Q)$:
\begin{longlist}[(2)]
\item[(1)] \textit{Regular phylogenies}:
$T_\theta\in\phy_{f,g}, \forall\theta\in[\Theta]$.

\item[(2)] \textit{Minimum frequency}:
$\nu_\theta\geq\underline{\nu}, \forall\theta\in[\Theta]$.
\end{longlist}
We denote by $\gmassump[f, g, \minnu, n]$
the set of $\Theta$-mixture models on $n$ leaves
satisfying these conditions.
\end{assumption}
%
\begin{remark}[(No minimum frequency)]
See the concluding remarks
for an extension of our techniques
when the minimum frequency assumption is
not satisfied.
\end{remark}

\subsubsection*{Tree identifiability}
Our first result states that,
under Assumption~\ref{assumpgeneral},
$\Theta$-mixture models are identifiable---except
for an ``asymptotically negligible fraction.''
To formalize this notion,
we use the following definition.
Note that $\gmassump[f,g,\minnu,n]$
is a compact subset of a finite product of metric spaces~\cite{BiHoVo01}
which we equip with its Borel $\sigma$-algebra.
%
\begin{definition}[(Permutation-invariant measure)]
Let
\[
A \subseteq\gmassump(f,g,\minnu,n)
\]
be a Borel set.
Given $\Theta$ permutations $\Pi= \{\Pi_\theta\}_{\theta\in
[\Theta]}$
of $X$, we let
\[
\Pi[\mixed]
\equiv\{\Pi_\theta[T_\theta]\}_{\theta\in[\Theta]}
\equiv\{(V_\theta,E_\theta;\phi_\theta\circ\Pi_\theta;\weight
_\theta
)\}_{\theta\in[\Theta]},\vadjust{\goodbreak}
\]
where $\circ$ indicates composition, and
\[
A_\Pi
= \{(\mixed,\nu,Q) \in\gmassump(f,g,\minnu,n)
\dvtx(\Pi[\mixed],\nu,Q) \in A\}.
\]
A probability measure $\lambda$ on
$\gmassump(f,g,\minnu,n)$ is
\textit{permutation-invariant} if for
all $A$ and $\Pi$ as above,
we have the following:
\[
\lambda[A]
= \lambda[A_\Pi].
\]
\end{definition}
%
\begin{remark}
Alternatively one can think of a
permutation-invariant measure as
first picking unlabeled trees, branch weights
and mixture frequencies according
to a specified joint distribution, and then
labeling the leaves of each tree in the
mixture independently, uniformly
at random. Note that the \textit{independent}
labeling of the trees is needed for our proof.
It ensures that the phylogenies in the
mixture are typically, ``sufficiently
distinct.'' Generalizing our results, possibly
in a weaker form,
to mixtures of ``similar'' phylogenies is
an important open problem. See~\cite{MosselRoch11a}
for recent progress in this direction.
\end{remark}

For two $\Theta$-mixture models
$(\mixed,\nu,Q)$ and
$(\mixed' = \{T'_\theta\}_{\theta\in[\Theta]},\nu',Q)$,
we write
\[
(\mixed,\nu,Q) \nsim(\mixed',\nu',Q),
\]
if there is no bijective mapping $h$
of $[\Theta]$ such that
\[
\tcal_l[T_{\theta}] =
\tcal_l\bigl[T'_{h(\theta)}\bigr]\qquad \forall\theta\in[\Theta].
\]
In words,
$(\mixed,\nu,Q)$ and $(\mixed',\nu',Q)$
are \textit{not} equivalent
up to component re-labeling.
%
\begin{theorem}[(Tree identifiability)]
\label{thm1}
Fix
$0 < f \leq g < +\infty$,
and $\minnu> 0$.
Then, there exists
a sequence of Borel
subsets
\[
A_n \subseteq\gmassump(f,g,\minnu,n),\qquad n \geq1,
\]
such that the following hold:
\begin{longlist}[(2)]
\item[(1)] For any sequence of
permutation-invariant measures
$\lambda_n$, $n \geq1$,
respectively, on $\gmassump(f,g,\minnu,n)$,
$n \geq1$, we have
\[
\lambda_n[A_n] = 1 - o_n(\minnu,f,g)
\]
as $n \to\infty$.
Here $o_n(\minnu,f,g)$ indicates
convergence to $0$ as $n \to\infty$
for fixed $\minnu,f,g$.

\item[(2)]
For all
\[
(\mixed,\nu,Q) \nsim(\mixed',\nu',Q)
\in\bigcup_{n \geq1} A_n,
\]
we have
\[
\dcal_l[(\mixed,\nu,Q)]
\neq\dcal_l[(\mixed',\nu',Q)].
\]
\end{longlist}
\end{theorem}
%
\begin{remark}
As remarked above, our proof requires that the phylogenies in the
mixture are ``sufficiently different.'' This is typically the case
under a permutation-invariant measure. Roughly speaking, the
complements of the sets $A_n$ in the previous theorem contain those
exceptional instances where the phylogenies are too ``similar.'' See
the proof for a formal definition of $A_n$.
\end{remark}

\subsubsection*{Tree distance}
We also generalize to GTR models a result of
Steel and Sz\'ekely:
phylogenies are typically far away
in variational distance~\cite{SteelSzekely06}. The techniques
in~\cite{SteelSzekely06}
apply only to group-based models and other
highly symmetric models; see~\cite{SteelSzekely06} for details.
Let \mbox{$\tv{\cdot}$} denote total variation distance; that is,
for two probability measures $\dcal$, $\dcal'$
on a measure space $(\Omega, \fcal)$ define
\[
\tv{\dcal- \dcal'}
= \sup_{B \in\fcal} |
\dcal(B) - \dcal'(B)
|.
\]

\begin{theorem}[(Tree distance)]\label{thm2}
Let $\{A_n\}_n$ be as in Theorem~\ref{thm1}
where $\Theta= 2$ and $\minnu= 1/2$
[in which case we necessarily have $\nu= (1/2,1/2)$]. Then for all
\[
(\mixed,\nu,Q) \in\bigcup_{n \geq1} A_n,
\]
we have
\[
\tv{\dcal_l[T_1,Q] - \dcal_l[T_2,Q]}
= 1 - o_n(1).
\]
\end{theorem}
%
\begin{remark}
Note that $\nu$ plays no substantive role in the previous theorem other
than to determine $A_n$.
\end{remark}

\subsubsection*{Tree reconstruction} The proof
of Theorems~\ref{thm1} and~\ref{thm2}
rely on the following reconstruction result
of independent interest.
We show that
the topologies can be reconstructed efficiently
with high confidence
using polynomial length sequences.
Recall that $k$ denotes the sequence length.
%
\begin{theorem}[(Tree reconstruction)]\label{thm3}
Fix
$0 < f \leq g < +\infty$,
and $\minnu> 0$.
Then, there exists
a sequence of Borel
subsets
\[
A_n \subseteq\gmassump(f,g,\minnu,n),\qquad n \geq1,
\]
such that the following hold:
\begin{longlist}[(2)]
\item[(1)] For any sequence of
permutation-invariant measures
$\lambda_n$, $n \geq1$,
respectively, on $\gmassump(f,g,\minnu,n)$,
$n \geq1$, we have
\[
\lambda_n[A_n] = 1 - o_n(\minnu,f,g)
\]
as $n \to\infty$.

\item[(2)] For all
\[
(\mixed,\nu,Q) \in\bigcup_{n \geq1} A_n,\vadjust{\goodbreak}
\]
the topologies of $(\mixed,\nu,Q)$
can be reconstructed
in time polynomial in $n$ and $k$
using polynomially many
samples (i.e., $k$ is polynomial
in $n$) with probability
$1 - o_n(\minnu,f,g)$ under
the samples
and the randomness of the algorithm.
\end{longlist}
\end{theorem}
%
\begin{remark}
The subsets $\{A_n\}_n$ in Theorems~\ref{thm1}
and~\ref{thm3} are in fact the same.
\end{remark}

The rest of the paper is devoted to the proof of
Theorem~\ref{thm3} which implies
Theorems~\ref{thm1} and~\ref{thm2}.

\subsection{Proof overview}
\label{sectionoverview}

The proof of Theorem~\ref{thm3} relies
on the construction of a \textit{clustering statistic} that
discriminates between distinct phylogenies.

\subsubsection*{Clustering statistic}
Fix $0 < f \leq g < +\infty$ and $\minnu> 0$.
Suppose for now that $\Theta= 2$, and
let $\lambda$ be a permutation-invariant
probability measure on $\gmassump[f,g$, $\minnu, n]$.
It will be useful to think of $\lambda$ as
a two-step procedure: first pick
unlabeled, weighted topologies; and
second, assign a uniformly random labeling
to the leaves of each tree.
Pick a $\Theta$-mixture model $(\mixed,\nu,Q)$
according to $\lambda$.
We will denote by $\P_\l$ and $\E_\l$
probability and expectation under~$\lambda$.
Similarly, we denote by $\P_l$ and $\E_l$
(resp., $\P_\A$ and $\E_\A$) probability
and expectation under $(\mixed,\nu,Q)$
(resp., under the randomness of our algorithm),
as well as combinations such as $\P_{\A,\l}$
with the obvious meaning.

Let $\mathbf{z} = (z_x)_{x=1}^r$ be a (real-valued) right
eigenvector of $Q$ corresponding to eigenvalue
$\Lambda_2 = -1$ and normalize $\mathbf{z}$ so that
\[
\sum_{x=1}^r \pi_x z_x^2 = 1.
\]
(Any negative eigenvalue could be used instead.)
Consider the following one-dimensional mapping
of the samples (\cite{MosselPeres03}, Lemma 5.3): for all $i =
1,\ldots,
k$ and
$a \in X$,
%
\begin{equation}\label{eqmapping}
\sigma_a^i = z_{s_a^i}.
\end{equation}
Recall that we drop the superscript when referring
to a single sample. It holds that
%
\begin{equation}\label{eqmeanzero}
\E_l[\sigma_a|N=\theta] = 0.
\end{equation}
Moreover, following a computation in~\cite{MosselPeres03}, Lemma 5.3,
letting $a \land b$ be the most recent common
ancestor of $a$ and $b$ (under the arbitrary
choice of root~$\rho$) one has
%
\begin{eqnarray} \label{eqcorr1}
q_\theta(a,b)
&=& \E_l[\sigma_a \sigma_b | N=\theta]
- \E_l[\sigma_a|N=\theta] \E[\sigma_b|N=\theta] \nonumber\\
&=& \E_l[\sigma_a \sigma_b | N=\theta]\nonumber\\
&=& \sum_{x=1}^r \pi_x
\E_l[\sigma_a \sigma_b | N=\theta, s_{a\land b} =
x]\nonumber\\[-8pt]\\[-8pt]
&=& \sum_{x=1}^r \pi_x
\E_l[\sigma_a | N=\theta, s_{a\land b} = x]\E_l[\sigma_b |
N=\theta,
s_{a\land b} = x]
\nonumber\\
&=& \sum_{x=1}^r \pi_x
\bigl(e^{-d_\theta(a\land b,a)}z_x\bigr)
\bigl(e^{-d_\theta(a\land b,b)}z_x\bigr)
\nonumber\\
&=& e^{-d_\theta(a,b)}\nonumber
\end{eqnarray}
and
%
\begin{equation}\label{eqcorr2}
q(a,b)
= \E_l[\sigma_a \sigma_b] - \E_l[\sigma_a] \E_l[\sigma_b]
= \E_l[\sigma_a \sigma_b]
= \sum_{\theta=1}^\Theta\nu_\theta e^{-d_\theta(a,b)}.
\end{equation}
We use a statistic of the form
%
\begin{equation}\label{eqclustering}
\ucal= \frac{1}{|\Upsilon|}\sum_{(a,b) \in\Upsilon} \sigma_a
\sigma_b,
\end{equation}
where $\Upsilon\subseteq\nsqneq$.
For $\ucal$ to be effective in discriminating
between $T_1$ and $T_2$, we require the following
(informal) conditions:
\begin{longlist}[(C3)]
\item[(C1)] The difference in conditional expectations
\[
\Delta= \bigl|\E_l[\ucal|N=1]
- \E_l[\ucal|N=2]\bigr|
\]
is large.

\item[(C2)] The statistic $\ucal$ is concentrated
around its mean
under both $\dcal_l[T_1,Q]$ and
$\dcal_l[T_2,Q]$.

\item[(C3)] The set $\Upsilon$ can be constructed
from data generated by the mixture $(\mixed,\nu,Q)$.
\end{longlist}
A $\ucal$ satisfying C1--C3 could be used to
infer the \textit{hidden} variables $N^1,\ldots,N^k$
and, thereby, to cluster the samples in their
respective component.

\subsubsection*{Prior work}
In~\cite{MosselRoch11a}, it was shown
in a related context that
taking $\Upsilon= \nsqneq$
is not in general an appropriate choice, as it may lead to a
large variance. Instead, the following lemma
was used.
%
\begin{uclaim*}[(Disjoint close pairs~\cite{SteelSzekely06}; see
also~\cite{MosselRoch11a})]
For any $T \in\phy_{f,g}^{(n)}$, there exists
a subset $\Upsilon\subseteq\nsqneq$ such that
the following hold:
\begin{longlist}[(2)]
\item[(1)] $|\Upsilon| = \Omega(n)$;

\item[(2)] $\forall(a,b) \in\Upsilon$,
$d_T(a,b) \leq3g$;

\item[(3)] $\forall(a_1,b_1)\neq(a_2,b_2) \in\Upsilon$,
the paths $\path_T(a_1,b_1)$ and $\path_T(a_2,b_2)$
are edge-disjoint. We will say that such pairs
are \textit{$T$-disjoint}.\vadjust{\goodbreak}
\end{longlist}
\end{uclaim*}

For special $Q$ matrices,
it was shown in~\cite{SteelSzekely06}
and~\cite{MosselRoch11a} that
such a $\Upsilon$ for $T = T_1$, say,
can be used to construct a
clustering statistic [similar to (\ref{eqclustering})]
concentrated under $\dcal_l[T_1,Q]$.
In particular, the $T_1$-disjointness assumption above
implies the independence of the variables
$\sigma_{a_1}\sigma_{b_1}$ and
$\sigma_{a_2}\sigma_{b_2}$ under
the $Q$ matrices considered
in~\cite{SteelSzekely06,MosselRoch11a}. Moreover,
Steel and Sz\'ekely~\cite{SteelSzekely06}
proved the existence of a further subset
that is also $T_2$-disjoint, but their construction
requires the knowledge of $T_2$.
Here we show how to satisfy conditions C1--C3
under GTR models.

\subsubsection*{High-level construction}
We give a sketch of our techniques.
Formal statements and full proofs can be
found in Sections~\ref{sectionpreliminaries},
\ref{sectionconstructing} and~\ref{sectiontree}.
For $\alpha> 0$, let
\[
\Upsilon_{\alpha, \theta}
= \{(a,b) \in\nsqneq\dvtx\dist_\theta(a,b) \leq\alpha\}
\]
and
\[
\Upsilon_{\alpha}
= \bigcup_{\theta\in[\Theta]} \Upsilon_{\alpha,\theta}.
\]
Because the variables $N^1, \ldots, N^k$ are hidden,
we cannot infer $\Upsilon_{\alpha,\theta}$ directly
from the samples, for instance,
using (\ref{eqcorr1}). Instead:

\mbox{}

\begin{center}
\fbox{
{\fontsize{9pt}{11pt}\selectfont{\begin{tabular}{p{325pt}}
(\textit{Step} 1)
Using (\ref{eqcorr2}) and the estimator
\[
\hat q(a,b) =
\frac{1}{k} \sum_{i=1}^k \sigma^i_a \sigma^i_b,
\]
we construct a set with size linear in $n$ satisfying
\[
\Upsilon_{4g} \subseteq\Upsilon'
\subseteq\Upsilon_{C_c}
\]
for an appropriate constant $C_c$; see Lemma~\ref{lemmaquasicherries}.
\end{tabular}}}
}
\end{center}

\mbox{}

Define
\[
\Upsilon'_\theta= \Upsilon' \cap\Upsilon_{C_c,\theta}.
\]
For general GTR rate matrices, $T_\theta$-disjointness
of $(a_1,b_1), (a_2,b_2) \in\Upsilon'_\theta$
does not guarantee independence
of $\sigma_{a_1}\sigma_{b_1}$ and
$\sigma_{a_2}\sigma_{b_2}$
under $\dcal_l[T_\theta,Q]$. Instead,
we choose pairs that are far enough from each other
by picking
a sufficiently sparse random subset of $\Upsilon'$; see Lemma \ref
{lemmaindependence}.
We say that $(a_1,b_1), (a_2,b_2) \in\Upsilon'_\theta$
are \textit{$T_\theta$-far} if the smallest
evolutionary distance between
$\{a_1,b_1\}$ and $\{a_2,b_2\}$
is at least $C_f \log\log n$
for a constant $C_f > 0$ to be determined.

\mbox{}

\begin{center}
\fbox{
{\fontsize{9pt}{11pt}\selectfont{\begin{tabular}{p{325pt}}
(\textit{Step} 2)
We take a random subset
$\Upsilon''$ of $\Upsilon'$ with
\[
|\Upsilon''| = \Theta(\log n);
\]
see Lemma~\ref{lemmasparsification2}.
\end{tabular}}}
}
\end{center}


Denoting
\[
\Upsilon''_\theta=
\Upsilon'' \cap\Upsilon_{C_c,\theta},
\]
we show that all
$(a_1,b_1) \neq(a_2,b_2) \in\Upsilon''_\theta$
are $T_\theta$-far.
Under a permutation-invariant~$\lambda$,
a pair $(a,b) \in\Upsilon_{\alpha,1}$ is
unlikely to be in $\Upsilon_{\alpha,2}$.
In particular,
we show that, under $\lambda$, the intersection of
$\Upsilon''_1$ and $\Upsilon''_2$ is empty.
In fact, a pair $(a,b) \in\Upsilon_{\alpha,1}$
is likely to be such that $d_2(a,b)$ is large.
We say that $(a,b) \in\nsqneq$ is
\textit{$T_\theta$-stretched} if
$d_\theta(a,b) \geq C_{st} \log\log n$
for a constant $C_{st} > 0$ to be determined.
We show that all $(a,b) \in\Upsilon''_1$
are $T_2$-stretched; see Lemma~\ref{lemmasparsification}.

To infer $\Upsilon''_\theta$, we consider
the quantity
\[
\hat r(c_1,c_2)
= \frac{1}{k} \sum_{i=1}^k
[
\sigma^i_{a_1} \sigma^i_{b_1} \sigma^i_{a_2} \sigma^i_{b_2}
- \hat q(a_1,b_1) \hat q(a_2,b_2)
]
\]
for $c_1 = (a_1,b_1) \neq c_2 = (a_2,b_2) \in\nsqneq$.
We note that if $(a,b) \in\Upsilon''$ is $T_2$-stretched,
then
\[
\E_l[\sigma_a \sigma_b|N=2] \approx\E_l[\sigma_a|N=2]
\E_l[\sigma_b|N=2] = 0
\]
and
\[
q(a,b) \approx\nu_1 q_1(a,b).
\]
There are then two cases:
\begin{longlist}[(II)]
\item[(I)] If
$c_1 = (a_1,b_1) \neq c_2 = (a_2,b_2) \in\Upsilon''_1$
(and similarly for $\Upsilon''_2$),
they are $T_1$-far and each is $T_2$-stretched.
Moreover
we show that $(c_1,c_2)$ is $T_2$-far.
Therefore,
\[
q(a_1,b_1) \approx\nu_1 q_1(a_1,b_1),
q(a_2,b_2) \approx\nu_1 q_1(a_2,b_2),
\]
and we show further that
\begin{eqnarray*}
&&
\E_l[\sigma_{a_1} \sigma_{b_1} \sigma_{a_2} \sigma_{b_2}]\\
&&\qquad\approx
\nu_1
\E_l[\sigma_{a_1} \sigma_{b_1}|N=1]
\E_l[\sigma_{a_2} \sigma_{b_2}|N=1]\\
&&\qquad\quad{} + \nu_2
\E_l[
\sigma_{a_1}|N=2] \E_l[\sigma_{b_1}|N=2]
\E_l[\sigma_{a_2}|N=2] \E_l[\sigma_{b_2}|N=2]\\
&&\qquad\approx \nu_1 q_1(a_1,b_1) q_1(a_2,b_2).
\end{eqnarray*}
So
\[
\hat r(c_1,c_2)
\approx\nu_1(1-\nu_1) q_1(a_1,b_1) q_1(a_2,b_2)
> 0.
\]

\item[(II)] On the other hand,
if $c_1 = (a_1,b_1) \in\Upsilon''_1$ and
$c_2 = (a_2,b_2) \in\Upsilon''_2$,
then $c_1$ is $T_2$-stretched, and
$c_2$ is $T_1$-stretched. Moreover we show that
$(c_1,c_2)$ is both $T_1$-far and
$T_2$-far. Therefore,
\[
q(a_1,b_1) \approx\nu_1 q_1(a_1,b_1),\qquad
q(a_2,b_2) \approx\nu_2 q_2(a_2,b_2),
\]
and we show that
\begin{eqnarray*}
\E_l[
\sigma_{a_1} \sigma_{b_1} \sigma_{a_2} \sigma_{b_2}
]
&\approx&
\nu_1 \E_l[
\sigma_{a_1} \sigma_{b_1}|N=1]
\E_l[\sigma_{a_2}|N=1] \E_l[\sigma_{b_2}|N=1]\\
&&{} + \nu_2 \E_l[
\sigma_{a_1}|N=2] \E_l[\sigma_{b_1}|N=2]
\E_l[\sigma_{a_2} \sigma_{b_2}|N=2]\\
&\approx& 0.
\end{eqnarray*}
So
\[
\hat r(c_1,c_2)
\approx- \nu_1 q_1(a_1,b_1) \nu_2 q_2(a_2,b_2) < 0;
\]
see Lemma~\ref{lemmaexpectations}.
\end{longlist}

The argument above leads to the following step.

\mbox{}

\begin{center}
\fbox{
{\fontsize{9pt}{11pt}\selectfont{\begin{tabular}{p{325pt}}
(\textit{Step} 3)
For all pairs $c_1 = (a_1,b_1)$
and $c_2 = (a_2,b_2)$ in $\Upsilon''$,
we compute $\hat r(c_1,c_2)$.
Using cases I and II, we then infer the sets
$\Upsilon''_1$ and
$\Upsilon''_2$. We form
the clustering statistics
\[
\ucal_\theta^i
=
\frac{1}{|\Upsilon''_\theta|}
\sum_{(a,b) \in\Upsilon''_\theta}
\sigma_a^i \sigma_b^i
\]
for $\theta= 1,2$ and
$i = 1,\ldots,k$; see Lemma~\ref{lemmaclusters}.
\end{tabular}}}
}
\end{center}

\mbox{}

By the arguments in cases I and II above,
we get that for $(a,b) \in\Upsilon''_1$,
\[
\E_l[\sigma_a \sigma_b|N=1]
\approx\nu_1 q_1(a,b),
\]
whereas
\[
\E_l[\sigma_a \sigma_b|N=2]
\approx\E_l[\sigma_a|N=2] \E_l[\sigma_b|N=2]
\approx0,
\]
so that (dropping the superscript
to refer to a single sample)
\[
\E_l[\ucal_1|N=1]
> C_\Delta,
\]
whereas
\[
\E_l[\ucal_1|N=2]
< C_\Delta
\]
for a constant $C_\Delta> 0$ to be determined
later; see Lemma~\ref{lemmaseparation}.
Moreover, the properties of $\Upsilon''_\theta$
discussed in cases I and II allow us to prove further
that $\ucal_\theta$ is concentrated around its
mean; see Lemma~\ref{lemmaconcentration}.
This leads to the following step.

\mbox{}

\begin{center}
\fbox{
{\fontsize{9pt}{11pt}\selectfont{\begin{tabular}{p{325pt}}
(\textit{Step} 4)
Divide the samples $i=1,\ldots,k$
into two clusters $K_1$ and $K_2$, according to whether
\[
\ucal^i_{1} > C_\Delta\quad\mbox{or}\quad
\ucal^i_{2} > C_\Delta,
\]
respectively; see Lemma~\ref{lemmabinning}.
\end{tabular}}}
}
\end{center}


Once the samples are divided into pure components,
we apply standard reconstruction
techniques to infer each topology.

\mbox{}

\begin{center}
\fbox{
{\fontsize{9pt}{11pt}\selectfont{\begin{tabular}{p{325pt}}
(\textit{Step} 5)
For $\theta=1,2$, reconstruct
the topology $\tcal_l[T_\theta]$
from the samples in $K_\theta$; see Lem\-ma~\ref{lemmadistortion}.
\end{tabular}}}
}
\end{center}

\subsubsection*{General $\Theta$} When $\Theta> 2$,
we proceed as above and construct a clustering
statistic for each component.

\section{Main lemmas}
\label{sectionpreliminaries}

In this section, we derive a number of
preliminary results. These results are also
described informally in Section~\ref{sectionoverview}.

Fix a GTR matrix $Q$ and constants
$\Theta\geq2$,
$0 < f \leq g < +\infty$ and $\minnu> 0$.
Let $\lambda$ be a permutation-invariant
probability measure on $\gmassump[f,g,\minnu,n]$.
Pick a $\Theta$-mixture model $(\mixed,\nu,Q)$
according to $\lambda$, and generate
$k$ independent samples $\{s_X^i\}_{i=1}^k$ from
$\dcal_l[(\mixed,\nu,Q)]$.
We work with the mapping $\{\sigma_X\}_{i=1}^k$
defined in (\ref{eqmapping}).

Throughout we assume that
the number of samples is
$k = n^{C_k}$ for some
$C_k > 0$ to be fixed later.

\subsection{Useful lemmas}

We will need the following standard concentration inequalities; see,
for example,~\cite{MotwaniRaghavan95}:
%
\begin{lemma}[(Azuma--Hoeffding inequality)]\label{lemmaazuma}
Suppose $\mathbf{Z}=(Z_1,\ldots,Z_m)$
are independent random variables taking values in a set
$S$, and $h\dvtx S^m \to\real$ is any $t$-Lipschitz function:
$|h(\bfz) - h(\bfz')|\leq t$ whenever $\bfz, \bfz' \in S^m$
differ at just one coordinate. Then,
$\forall\zeta> 0$,
\[
\prob[|h(\mathbf{Z}) - \expec[h(\mathbf{Z})]| \geq\zeta]
\leq2\exp\biggl(-\frac{\zeta^2 }{2 t^2 m}\biggr).
\]
\end{lemma}
%
\begin{lemma}[(Chernoff bounds)]\label{lemmachernoff}
Let $Z_1,\ldots,Z_m$ be independent Poisson trials
such that, for $1 \leq i \leq m$, $\P[Z_i = 1] = p_i$
where $0 < p_i < 1$. Then, for
$Z = \sum_{i=1}^m Z_i$, $M = \E[Z] = \sum_{i=1}^m p_i$,
$0 < \delta_- \leq1$, and $\delta_+ > 2e - 1$,
\[
\P[Z < (1-\delta_-) M] < e^{-M \delta_-^2/2}
\]
and
\[
\P[Z > (1+\delta_+) M] < 2^{-(1+\delta_+)M}.
\]
\end{lemma}

\subsection{Large-sample asymptotics}

Denoting $\kcal= [k]$, let
$\kcal_\theta\subseteq\kcal$ be those
samples coming from component $\theta$,
that is,
\[
\kcal_\theta= \{i\in\kcal\dvtx N^i = \theta\}.
\]

\begin{lemma}[(Size of $\kcal_\theta$)]
\label{lemmakcaltheta}
Under $\P_l$, for any $C_s > 1$, we have
\[
C_s^{-1}
\leq\frac{|\kcal_\theta|}{\nu_\theta k}
\leq C_s
\]
for all $\theta\in[\Theta]$,
except with probability $\exp(-\Omega(n^{C_k}))$.
\end{lemma}
\begin{pf}
Recall that $\minnu\leq\nu_\theta\leq1 - \minnu$.
Using Lemma~\ref{lemmaazuma} with
$m = k$ and
\[
\zeta= \nu_\theta k
\max\{1 - C_s^{-1}, C_s - 1\}
= \nu_\theta k (C_s - 1)
\]
gives the result.
\end{pf}

Consider the estimators
\[
\hat q_\theta(a,b) =
\frac{1}{|\kcal_\theta|}
\sum_{i\in\kcal_\theta} \sigma^i_a \sigma^i_b
\]
and
\[
\hat q(a,b) =
\frac{1}{k} \sum_{i=1}^k \sigma^i_a \sigma^i_b.
\]
Let
\[
q_\theta(a,b) = e^{-d_\theta(a,b)}
\]
and
\[
q(a,b) = \sum_{\theta=1}^\Theta\nu_\theta q_\theta(a,b).
\]

\begin{lemma}[(Accuracy of $\hat q$)]
\label{lemmahatq}
Fix $0 < C_q < C_k/2$.
Under $\P_l$, we have
\[
|\hat q(a,b) - q(a,b)|
\leq n^{-C_q}
\]
and
\[
|\hat q_\theta(a,b) - q_\theta(a,b)|
\leq n^{-C_q}
\]
for all $\theta\in[\Theta]$
and all $(a,b) \in\nsqneq$
except with probability
$\exp(-{\poly(n)})$.
\end{lemma}
\begin{pf}
For each $(a,b) \in\nsqneq$, $\hat q(a,b)$
is a sum of $k$ independent variables.
By Lemma~\ref{lemmaazuma}, taking
$m = k$,
$t = k^{-1} \max_i |z_i|^2$,
$\zeta= n^{-C_q}$,
we have
\[
|\hat q(a,b) - q(a,b)| \leq n^{-C_q},
\]
except\vspace*{1pt} with probability $2\exp(-\Omega(n^{C_k - 2 C_q}))$.
Note that there are at most $n^2$
elements in $\nsqneq$ so that
the probability of failure is at most
\[
2 n^2 \exp(- \Omega(n^{C_k - 2 C_q}))
= \exp(-\Omega(n^{C_k - 2 C_q})).
\]

Using Lemma~\ref{lemmakcaltheta},
the same holds for each $\theta$. The overall
probability of failure under $\P_l$ is
$\exp(-\Omega(n^{C_k - 2 C_q}))$.
\end{pf}

Following the same argument,
a similar result holds for
\[
\hat r(c_1,c_2)
= \frac{1}{k} \sum_{i=1}^k
[
\sigma^i_{a_1} \sigma^i_{b_1} \sigma^i_{a_2} \sigma^i_{b_2}
- \hat q(a_1,b_1) \hat q(a_2,b_2)
]
\]
for $c_1 = (a_1,b_1) \neq c_2 = (a_2,b_2) \in\nsqneq$.
Let
\[
r(c_1,c_2) = \E_l[\hat r(c_1,c_2)].
\]

\begin{lemma}[(Accuracy of $\hat r$)]
\label{lemmahatr}
Under $\P_l$, we have
\[
|\hat r(c_1,c_2) - r(c_1,c_2)|
\leq n^{-C_q}
\]
for all $c_1 = (a_1,
b_1) \neq c_2 = (a_2,b_2) \in\nsqneq$
except with probability
$\exp(-{\poly(n)})$.
\end{lemma}

\subsection{Combinatorial properties}
\label{seccombinatorial}

For $\alpha> 0$, let
%
\begin{equation}\label{equpsilontheta}
\Upsilon_{\alpha, \theta}
= \{(a,b) \in\nsqneq\dvtx\dist_\theta(a,b) \leq\alpha\}
\end{equation}
and
%
\begin{equation}\label{equpsilonalpha}
\Upsilon_{\alpha}
= \bigcup_{\theta\in[\Theta]} \Upsilon_{\alpha,\theta}.
\end{equation}
The lower bound below
follows from a (stronger) lemma in~\cite{SteelSzekely06}; see
also~\cite{MosselRoch11a}.

\begin{lemma}[(Size of $\Upsilon_{\alpha,\theta}$)]\label{lemmasize}
For all $\alpha> 0$ and $\theta\in[\Theta]$,
\[
\tfrac{1}{4} n
\leq
|\Upsilon_{\alpha,\theta}|
\leq
2^{\lfloor{\alpha}/{f}\rfloor} n.
\]
In particular,
\[
\tfrac{1}{4} n
\leq
|\Upsilon_{\alpha}|
\leq
\Theta2^{\lfloor{\alpha}/{f}\rfloor} n.
\]
\end{lemma}
\begin{pf}
For $a \in X$ and $\alpha\geq4g$, let
\[
\bcal_\alpha(a) = \{v \in V\dvtx\dist_\theta(\phi_\theta(a),v)
\leq\alpha\}.
\]
Since $T_\theta$ is binary, there are at most
$2^{\lfloor{\alpha}/{f}\rfloor}$ vertices within
evolutionary distance~$\alpha$, that is,
\[
|\bcal_\alpha(a)| \leq
2^{\lfloor{\alpha}/{f}\rfloor}.
\]
Restricting to leaves gives the upper bound.

Let
\[
\Gamma_\alpha= \bigl\{a \in[n]\dvtx\dist_\theta(a,b) > \alpha,
\forall
b\in[n]-\{a\}\bigr\},
\]
that is, $\Gamma_\alpha$ is the set of leaves with no other leaf at
evolutionary distance $\alpha$ in~$T_\theta$.
We will bound the size of $\Gamma_\alpha$.
Note that for all $a,b \in\Gamma_\alpha$ with $a\neq b$, we have
$\bcal_{{\alpha}/{2}}(a) \cap\bcal{{\alpha}/{2}}(b) =
\varnothing$
by the triangle inequality. Moreover, it holds that for all \mbox{$a \in
\Gamma_\alpha$}
\[
|\bcal_{{\alpha}/{2}}(a)| \geq2^{\lfloor{\alpha}/({2g})\rfloor},
\]
since $T_\theta$ is binary, and there is no leaf other than $a$ in
$\bcal_{{\alpha}/{2}}(a)$.
Hence, we must have
\[
|\Gamma_\alpha| \leq\frac{2n - 2}{2^{\lfloor{\alpha}/({2g})\rfloor}}
\leq\biggl(\frac{1}{2^{\lfloor{\alpha}/({2g})
\rfloor-
1}}\biggr) n
\]
as there are $2n - 2$ nodes in $T_\theta$.
Now, for all $a \notin\Gamma_\alpha$ assign an arbitrary leaf at
evolutionary distance at most $\alpha$.
Then
\begin{eqnarray*}
|\Upsilon_{\alpha,\theta}|
&\geq& \frac{1}{2}(n - |\Gamma_\alpha|)\\
&\geq& \frac{1}{2}\biggl(1 - \frac{1}{2^{\lfloor{\alpha}/({2g})\rfloor- 1}}\biggr) n,
\end{eqnarray*}
where we divided by $2$ to avoid double-counting.
The result follows from the assumption $\alpha\geq4g$.
\end{pf}

Let $C_c > 4g$, $C_f > 0$, and $C_{st} > C_f$
to be fixed later.
%
\begin{definition}[($T_\theta$-quasicherry)]
We say that $(a,b) \in\nsqneq$ is a
\textit{$T_\theta$-quasi\-cherry}
if $(a,b) \in\Upsilon_{C_c,\theta}$.
\end{definition}
%
\begin{definition}[($T_\theta$-stretched)]
We say that $(a,b) \in\nsqneq$ is
\textit{$T_\theta$-stretched} if
$d_\theta(a,b) \geq C_{st} \log\log n$.
\end{definition}
%
\begin{definition}[($T_\theta$-far)]
We say that $c_1 = (a_1,b_1) \neq c_2 = (a_2,b_2) \in\nsqneq$
are \textit{$T_\theta$-far} if
\[
d_\theta(c_1,c_2) \equiv\min\bigl\{d_\theta(x_1,x_2)\dvtx x_1 \in\{a_1,b_1\},
x_2 \in\{a_2,b_2\}\bigr\} \geq C_f \log\log n.
\]
\end{definition}

Let $\Upsilon'$ be
any subset satisfying
%
\begin{equation}\label{eqbetween}
\Upsilon_{4g} \subseteq\Upsilon' \subseteq\Upsilon_{C_c}
\end{equation}
and let
%
\begin{equation}\label{eqrestriction}
\Upsilon'_\theta= \Upsilon' \cap\Upsilon_{C_c,\theta}.
\end{equation}
Let $C_{sp}^p > 0$ to be fixed later.
Keep
each $(a,b) \in\Upsilon_{C_c}$ independently
with probability
\[
p_{sp} = \frac{C_{sp}^p \log n}{n}
\]
to form the set $\Upsilon''_{C_c}$,
and let
\[
\Upsilon'' = \Upsilon' \cap\Upsilon''_{C_c}.
\]
Let $0 < C_{sp}^- < C_{sp}^+ < +\infty$
be constants (to be determined).
%
\begin{definition}[(Properly sparse)]
\label{defsparsification}
A subset
$\Upsilon_{4g} \subseteq\Upsilon'' \subseteq\Upsilon_{C_c}$
with
\[
\Upsilon''_\theta= \Upsilon'' \cap\Upsilon_{C_c,\theta},\qquad
\theta\in[\Theta],
\]
is \textit{properly sparse} if
it satisfies the following properties:
For all $\theta\in[\Theta]$:
\begin{longlist}[(2)]
\item[(1)] We have
$
C_{sp}^- \log n
\leq|\Upsilon''_\theta|
\leq C_{sp}^+ \log n.
$

\item[(2)] All
$c_1 = (a_1,b_1) \neq c_2 = (a_2,b_2) \in\Upsilon''$
are $T_\theta$-far.

\item[(3)] All
pairs in $\Upsilon''_\theta$
are $T_{\theta'}$-stretched for $\theta' \neq\theta$.
\end{longlist}
\end{definition}

Let
\[
\Upsilon''_{C_c,\theta} = \Upsilon''_{C_c} \cap\Upsilon
_{C_c,\theta},\qquad
\theta\in[\Theta],
\]
and
\[
\Upsilon''_{4g,\theta} = \Upsilon_{4g}
\cap\Upsilon''_{C_c,\theta},\qquad
\theta\in[\Theta].
\]

\begin{lemma}[(Sparsification)]
\label{lemmasparsification}
There exist constants
$0 < C_{sp}^- < C_{sp}^+ < +\infty$ such that,
under $\P_{\A,\l}$, the set $\Upsilon''_{C_c}$
as above
satisfies\vspace*{1pt} the following properties, except
with probability $1/\poly(n)$:
for all $\theta\in[\Theta]$:
\begin{longlist}[(2)]
\item[(1)] We have
$
C_{sp}^- \log n
\leq|\Upsilon''_{4g,\theta}|$
and
$
|\Upsilon''_{C_c,\theta}|
\leq C_{sp}^+ \log n.
$

\item[(2)] All
$c_1 = (a_1,b_1) \neq c_2 = (a_2,b_2) \in\Upsilon''_{C_c}$
are $T_\theta$-far.

\item[(3)] All
pairs in $\Upsilon''_{C_c,\theta}$
are $T_{\theta'}$-stretched for $\theta' \neq\theta$.

\end{longlist}
In particular, the set $\Upsilon''$ as above is
properly sparse.
Moreover, the claim holds for any $C_{sp}^- > 0$
by taking $C_{sp}^p > 0$ large enough.
\end{lemma}

Intuitively, part (2) follows from the sparsification
step whereas part (3) is a consequence
of the permutation-invariance of $\l$. We give
a formal proof next.
\begin{pf*}{Proof of Lemma~\ref{lemmasparsification}}
For part (1), we use Lemma~\ref{lemmachernoff}.
Take
\[
\frac{1}{4}
C_{sp}^p \log n
\leq M_{4g}
\equiv\frac{C_{sp}^p \log n}{n} |\Upsilon_{4g,\theta}|
\]
and
\[
M_{C_c}
\equiv\frac{C_{sp}^p \log n}{n} |\Upsilon_{C_c,\theta}|
\leq
2^{\lfloor{C_c}/{f}\rfloor} C_{sp}^p \log n.
\]
With
$\delta_- = 1/2$,
$\delta_+ = 5$,
we have
\[
\P_\A[|\Upsilon''_{4g,\theta}| <(1-\delta_-) M_{4g}] < e^{-M_{4g}
\delta_-^2/2}
= \frac{1}{\poly(n)}\vadjust{\goodbreak}
\]
and
\[
\P_\A[|\Upsilon''_{C_c,\theta}| > (1+\delta_+) M_{C_c}] <
2^{-(1+\delta
_+)M_{C_c}}
= \frac{1}{\poly(n)}.
\]
The first part follows from the choice
\[
C_{sp}^-
= \frac{C_{sp}^p}{8}
\]
and
\[
C_{sp}^+
= 6 C_{sp}^p 2^{\lfloor{C_c}/{f}\rfloor}.
\]

For the second part, let $c_1 = (a_1,b_1)$ be a pair
in $\Upsilon''_{C_c}$.
Let $\scal$ be the collection of pairs $c_2 = (a_2,b_2)
\neq c_1$
in the original set $\Upsilon_{C_c}$ that are
within evolutionary
distance $C_f \log\log n$ of $c_1$ in $T_\theta$,
that is,
\[
d(c_1,c_2) \leq C_f \log\log n.
\]
Note that the number of leaves
within evolutionary distance
$C_f \log\log n$ from $a_1$ or $b_1$
is at most $2 \cdot2^{\lfloor{C_f \log\log n}/{f}
\rfloor}$. Moreover,
each such leaf can be involved in at most
$\Theta2^{\lfloor{C_c}/{f}\rfloor}$ pairs,
since any pair in $\Upsilon_{C_c}$ must be a $T_{\theta'}$-quasicherry
for some $\theta' \in[\Theta]$
and the number of leaves at evolutionary
distance $C_c$ from a vertex in a tree in
$\phy_{f,g}$ is at most $2^{\lfloor{C_c}/{f}\rfloor
}$. Hence
\[
|\scal| \leq2 \cdot2^{\lfloor{C_f \log\log n}/{f}
\rfloor} \cdot\Theta2^{\lfloor{C_c}/{f}\rfloor}
= O(\log n).
\]
Therefore the probability that any $c_2 \in\scal$
remains in $\Upsilon''_{C_c}$ is at most $O(\log^2 n/\break n)$.
Assuming part (1) holds, summing over $\Upsilon''_{C_c}$,
and applying Markov's inequality, we get
\[
\P_\A[|c_1 \neq c_2 \in\Upsilon''_{C_c}\dvtx c_1, c_2\mbox{ are not
$T_\theta$-far} | \geq1]
= O\biggl(\frac{\log^3 n}{n}\biggr) + \frac{1}{\poly(n)}.
\]
This gives the second part.

For the third part, consider a
$T_\theta$-quasicherry $(a,b)$.
Thinking of $\l$ as assigning
leaf labels in $T_{\theta'}$ uniformly at random,
the probability that $b$ is within evolutionary
distance $C_{st} \log\log n$ of $a$ in $T_{\theta'}$
is at most
\[
\P_\l[\mbox{$(a,b)$ is not $T_{\theta'}$-stretched}]
\leq\frac{2^{\lfloor{C_{st} \log\log n}/{f}
\rfloor
}}{n} = O\biggl(\frac{\log n}{n}\biggr),
\]
where the numerator in the second expression
is an upper bound on the number of vertices
at evolutionary
distance $C_{st} \log\log n$ of $a$ in $T_{\theta'}$.
Summing over all pairs in $\Upsilon''_{C_c,\theta}$ and
assuming the bound in part (1) holds,
the expected number of pairs in $\Upsilon''_{C_c,\theta}$
that are not $T_{\theta'}$-stretched is
$O(\log^2 n/n)$. By Markov's inequality,
\[
\P_{\A,\l}[|\{
(a,b) \in\Upsilon''_{C_c,\theta}\dvtx (a,b)\mbox{ is not $T_{\theta
'}$-stretched}\}| \geq1] \leq
O\biggl(\frac{\log^2 n}{n}\biggr) + \frac{1}{\poly(n)}.
\]
This gives the third part.
\end{pf*}

\subsection{Mixing}

We use a mixing argument similar to~\cite{Mossel03}.
Let
\[
\qmin= \min_{x\neq y} Q_{x y},
\]
which is positive by assumption. We think of $Q$ as acting as follows.
From a state~$x$, we have two type of transitions to $y \neq x$:
\begin{longlist}[(ii)]
\item[(i)] We jump to state $y$ at rate $\qmin> 0$.

\item[(ii)] We jump to state $y$ at rate $Q_{x y} - \qmin\geq0$.
\end{longlist}
Note that a transition of type (i) does not depend on
the starting state. Hence if $\pcal$ is a path
from $u$ to $v$ in $T_\theta$, $N = \theta$,
and a transition of type (i) occurs along $\pcal$,
then $\sigma_u$ is \textit{independent} of $\sigma_v$.
The probability, conditioned on $N=\theta$, that such
a transition does not occur, is $e^{-d_\theta(u,v) (r-1) \qmin}$.

Let $\Upsilon'' \subseteq\nsqneq$
be a properly sparse set.
We show next that pairs in $\Upsilon''$
are independent with high probability.
We proceed by considering the paths
joining them and arguing that transitions
of type (i) are likely to occur on them
by the combinatorial properties
in Definition~\ref{defsparsification}.
Formally, fix $\theta\in[\Theta]$, and
consider two pairs $c_1 = (a_1,b_1)
\neq c_2 = (a_2,b_2) \in\Upsilon''$.
By Definition~\ref{defsparsification}, $c_1$ and
$c_2$ are $T_\theta$-far.
There are three cases without loss
of generality:
\begin{longlist}[(2)]
\item[(1)] \textit{$c_1, c_2$ are $T_\theta$-quasicherries.}
In the subtree
of $T_\theta$ connecting
$\{a_1,b_1,a_2,\break b_2\}$, called a \textit{quartet},
the paths $\path_{T_\theta}(a_1,b_1)$
and $\path_{T_\theta}(a_2,b_2)$ are disjoint.
This is denoted by the \textit{quartet split}
$a_1 b_1|a_2 b_2$. Let $\pcal^\theta[c_1,c_2]$
be the internal path of the quartet. Note that
by Definition~\ref{defsparsification}
the length of $\pcal^\theta[c_1,c_2]$ is
at least $C_f \log\log n - 2 C_c$.
Denote by $\pcal_{c_1}^\theta[c_1,c_2]$
the subpath of $\pcal^\theta[c_1,c_2]$
within evolutionary distance
$\frac{1}{3}C_f \log\log n$ of
$c_1$.

\item[(2)] \textit{$c_1$ is a $T_\theta$-quasicherry}, \textit{and
$c_2$ is $T_\theta$-stretched.}
Consider the subtree
of $T_\theta$ connecting
$\{a_1,b_1,a_2\}$, called a \textit{triplet}, and
let $u$ be the central vertex of it.
Let $\pcal^\theta[c_1,a_2]$
be the path connecting $u$ and $a_2$.
Note that
by Definition~\ref{defsparsification},
the length of $\pcal^\theta[c_1,a_2]$ is
at least $C_f \log\log n - C_c$.
Denote by $\pcal_{c_1}^\theta[c_1,a_2]$
the subpath of $\pcal^\theta[c_1,a_2]$
within evolutionary distance
$\frac{1}{3}C_f \log\log n$ of
$c_1$. Similarly,
denote by $\pcal_{a_2}^\theta[c_1,a_2]$
the subpath of $\pcal^\theta[c_1,a_2]$
within evolutionary distance
$\frac{1}{3}C_f \log\log n$ of
$a_2$.

\item[(3)] \textit{$c_1, c_2$ are $T_\theta$-stretched.}
Let $\pcal^\theta[a_1,a_2]$
be the path connecting $a_1$ and $a_2$.
Note that
by Definition~\ref{defsparsification}
the length of $\pcal^\theta[a_1,a_2]$ is
at least $C_f \log\log n$.
Denote by $\pcal_{a_1}^\theta[a_1,a_2]$
the subpath of $\pcal^\theta[a_1,a_2]$
within evolutionary distance
$\frac{1}{3}C_f \log\log n$ of
$a_1$. Similarly, let $\pcal^\theta[a_1,b_1]$
be the path joining $a_1$ and $b_1$, and
let $\pcal_{a_1}^\theta[a_1,b_1]$ be
the subpath of $\pcal^\theta[a_1,b_1]$
within evolutionary distance
$\frac{1}{3}C_{st} \log\log n
> \frac{1}{3}C_f \log\log n$ of
$a_1$.
\end{longlist}
Condition on $N = \theta$.
For each $c_1 = (a_1,b_1) \in\Upsilon''_\theta$,
let $\ecal_{c_1}^\theta$ be the following event:

\begin{quote}
Each subpath $\pcal_{c_1}^\theta[c_1,c_2]$,
$c_2\neq c_1 \in\Upsilon''_\theta$, and
each subpath $\pcal_{c_1}^\theta[c_1,a_2]$,
$c_2 = (a_2,b_2) \in\Upsilon'' - \Upsilon''_\theta$,
undergo a transition of type (i) during
the generation of sample $\sigma_X$.
\end{quote}

\noindent Similarly,
for each $c_1 = (a_1,b_1) \in
\Upsilon'' - \Upsilon''_\theta$,
let $\ecal_{c_1}^\theta= \ecal_{a_1}^\theta
\cap\ecal_{b_1}^\theta$
where $\ecal_{a_1}^\theta$ is the following event
(and similarly for $\ecal_{b_1}^\theta$):

\begin{quote}
Each subpath $\pcal_{a_1}^\theta[c_2,a_1]$,
$c_2 \in\Upsilon''_\theta$,
each subpath $\pcal_{a_1}^\theta[a_1,a_2]$,
$c_2 = (a_2,b_2) \in\Upsilon'' - \Upsilon''_\theta$
with $c_1 \neq c_2$,
as well as subpath $\pcal_{a_1}^\theta[a_1,b_1]$
undergo a transition of type (i) during
the generation of sample $\sigma_X$.
\end{quote}

\noindent Note that, under $\ecal_{c_1}^\theta$, the random
variable $\sigma_{a_1}\sigma_{b_1}$ is independent
of every other such random variable in $\Upsilon''$.
Moreover, in the case $c_1 \in\Upsilon'' - \Upsilon''_\theta$,
then further $\sigma_{a_1}$ is independent of
$\sigma_{b_1}$.
The next lemma shows that most of the events above
occur with high probability implying that a large
fraction of $\sigma_{a_1}\sigma_{b_1}$'s
are mutually independent.
%
\begin{lemma}[(Pair independence)]
\label{lemmaindependence}
Let $\Upsilon'' \subseteq\nsqneq$ be a
properly sparse set.
Conditioned on $N = \theta$,
let
\[
\ical
= \{c_1 \in\Upsilon''\dvtx\ecal_{c_1}^\theta\mbox{ holds}\}.
\]
For any $0 < \eps_\ical< 1$
and $C_\ical> 0$,
there exist $C_f$,
$C_{st} > C_f$ and $C_{sp}^- > 0$ large enough
so that
the following holds except
with probability $n^{-C_\ical}$
under $\P_l$:
\[
|\ical|
\geq(1 - \eps_\ical)|\Upsilon''|.
\]
\end{lemma}
\begin{pf}
Condition on $N = \theta$. Note that
the $\ecal_{c_1}^\theta$'s are mutually independent
because the corresponding paths are
disjoint by construction. By a union bound over
$\Upsilon''$,
for all $c_1 \in\Upsilon''$,
%
\begin{eqnarray} \label{eqball}
\P_l[(\ecal_{c_1}^\theta)^c| N = \theta]
&\leq& 2 C_{sp}^+ \log n \cdot e^{-(({1}/{3}) C_f \log\log n
- 2 C_c)
(r-1) \qmin}\nonumber\\[-8pt]\\[-8pt]
&=& \frac{1}{\poly(\log n)}\nonumber
\end{eqnarray}
for $C_f$ large enough. Applying
Lemma~\ref{lemmachernoff} with
\[
M = |\Upsilon''| \cdot\P_l[(\ecal_{c_1}^\theta
)^c| N =
\theta]
\]
and $\delta_+ > 2e$ such that
\[
(1+ \delta_+) M = \eps_\ical|\Upsilon''|
\geq\eps_\ical C_{sp}^- \log n,
\]
we get
\[
\P_l[|\Upsilon'' - \ical| > \eps_\ical|\Upsilon''|]
\leq2^{-\eps_\ical|\Upsilon''|}
= \frac{1}{n^{C_\ical}}
\]
by taking $C_{sp}^-$
large enough in Definition~\ref{defsparsification}.
\end{pf}

We use the independence claims above to
simplify expectation computations.
%
\begin{lemma}[(Expectation computations)]
\label{lemmaexpectations}
Let $\Upsilon'' \subseteq\nsqneq$ be a
properly sparse set.
The following hold. For all $\theta\neq\theta' \in[\Theta]$:
\begin{longlist}[(2)]
\item[(1)]
$\forall(a,b) \in\Upsilon''_\theta$,
\[
q_\theta(a,b)
\geq e^{-C_c}.
\]

\item[(2)]
$\forall(a,b) \in\Upsilon'' - \Upsilon''_\theta$,
\[
q_\theta(a,b)
= \frac{1}{\poly(\log n)}.
\]

\item[(3)]
$\forall(a,b) \in\Upsilon''_\theta$,
\[
q(a,b)
= \nu_\theta q_\theta(a,b)
+ \frac{1}{\poly(\log n)}.
\]

\item[(4)]
$\forall c_1 = (a_1,b_1)
\neq c_2 = (a_2,b_2) \in\Upsilon''_\theta$,
\begin{eqnarray*}
r(c_1,c_2)
&=& \nu_\theta(1 - \nu_\theta) q_\theta(a_1,b_1)
q_\theta(a_2,b_2) + \frac{1}{\poly(\log n)}\\
&\geq& \frac{1}{2} \minnu(1-\minnu) e^{-2C_c}
> 0.
\end{eqnarray*}

\item[(5)]
$\forall c_1 = (a_1,b_1) \in\Upsilon''_\theta,
c_2 = (a_2,b_2) \in\Upsilon''_{\theta'}$,
\begin{eqnarray*}
r(c_1,c_2)
&=& -\nu_\theta q_\theta(a_1,b_1)
\nu_{\theta'} q_{\theta'}(a_2,b_2) + \frac{1}{\poly(\log n)}\\
&\leq& -\frac{1}{2} \minnu e^{-2C_c}
< 0.
\end{eqnarray*}
\end{longlist}
\end{lemma}
\begin{pf}
Parts (1) and (2) follow
from the fact that $q_\theta(a,b)
= e^{-d_\theta(a,b)}$, $d_\theta(a,b) \leq C_c$
for all $(a,b) \in\Upsilon''_\theta$ and
$d_\theta(a,b) \geq C_{st} \log\log n$
for all $(a,b) \in\Upsilon'' - \Upsilon''_\theta$
from Definition~\ref{defsparsification}.
Part (3) follows from parts (1) and (2).

For part (4),
let
$c_1 = (a_1,b_1)
\neq c_2 = (a_2,b_2) \in\Upsilon''_\theta$.
Note that
\begin{eqnarray*}
\E_l[\sigma_{a_1}
\sigma_{b_1}
\sigma_{a_2}
\sigma_{b_2}
|
N = \theta,
\ecal_{c_1}^\theta,
\ecal_{c_2}^\theta]
&=&\E_l[\sigma_{a_1}
\sigma_{b_1}| N = \theta]
\E_l[\sigma_{a_2}
\sigma_{b_2}
|
N = \theta]\\
&=&
q_\theta(a_1,b_1)
q_\theta(a_2,b_2)
\end{eqnarray*}
and
\begin{eqnarray*}
\E_l[\sigma_{a_1}
\sigma_{b_1}
\sigma_{a_2}
\sigma_{b_2}
|
N = \theta',
\ecal_{c_1}^{\theta'},
\ecal_{c_2}^{\theta'}]
&=&
\E_l[\sigma_{a_1}|N = \theta']
\E_l[\sigma_{b_1}|N = \theta']\\
&&{} \times
\E_l[\sigma_{a_2}|N = \theta']
\E_l[\sigma_{b_2}|N = \theta']\\
&=& 0
\end{eqnarray*}
by (\ref{eqmeanzero}), so that
\[
\E_l[\sigma_{a_1}
\sigma_{b_1}
\sigma_{a_2}
\sigma_{b_2}]
= \nu_\theta q_\theta(a_1,b_1)
q_\theta(a_2,b_2)
+ \frac{1}{\poly(\log n)}
\]
from (\ref{eqball}).
Then part (4) follows from Lemma~\ref{lemmahatq}
and part (3).\vadjust{\goodbreak}

For part (5),
let
$c_1 = (a_1,b_1) \in\Upsilon''_\theta,
c_2 = (a_2,b_2) \in\Upsilon''_{\theta'}$.
Let $\theta'' \neq\theta, \theta'$.
Note that
\begin{eqnarray*}
\E_l[\sigma_{a_1}
\sigma_{b_1}
\sigma_{a_2}
\sigma_{b_2}
|
N = \theta,
\ecal_{c_1}^{\theta},
\ecal_{c_2}^{\theta}]
&=&
\E_l[\sigma_{a_1}
\sigma_{b_1}|N = \theta]\\
&&{} \times
\E_l[\sigma_{a_2}|N = \theta]
\E_l[\sigma_{b_2}|N = \theta]\\
&=& 0
\end{eqnarray*}
and
\begin{eqnarray*}
\E_l[\sigma_{a_1}
\sigma_{b_1}
\sigma_{a_2}
\sigma_{b_2}
|
N = \theta',
\ecal_{c_1}^{\theta'},
\ecal_{c_2}^{\theta'}]
&=&
\E_l[\sigma_{a_1}
\sigma_{b_1}|N = \theta']\\
&&{} \times
\E_l[\sigma_{a_2}|N = \theta']
\E_l[\sigma_{b_2}|N = \theta']\\
&=& 0.
\end{eqnarray*}
Moreover, since $c_1,c_2 \notin\Upsilon''_{\theta''}$,
\begin{eqnarray*}
\E_l[\sigma_{a_1}
\sigma_{b_1}
\sigma_{a_2}
\sigma_{b_2}
|
N = \theta'',
\ecal_{c_1}^{\theta''},
\ecal_{c_2}^{\theta''}]
&=&
\E_l[\sigma_{a_1}|N = \theta'']
\E_l[\sigma_{b_1}|N = \theta'']\\
&&{} \times
\E_l[\sigma_{a_2}|N = \theta'']
\E_l[\sigma_{b_2}|N = \theta'']\\
&=& 0.
\end{eqnarray*}
Hence
\[
\E_l[\sigma_{a_1}
\sigma_{b_1}
\sigma_{a_2}
\sigma_{b_2}]
= 0
+ \frac{1}{\poly(\log n)}
\]
from (\ref{eqball}).
Then part (5) follows from Lemma~\ref{lemmahatq}
and part (3).
\end{pf}

\subsection{Large-tree concentration}

Let $\Upsilon'' \subseteq\nsqneq$ be a
properly sparse set.
Consider the clustering statistic
\[
\ucal_\theta
=
\frac{1}{|\Upsilon''_\theta|}
\sum_{(a,b) \in\Upsilon''_\theta}
\sigma_a \sigma_b.
\]
We show that $\ucal_\theta$ is concentrated
and separates the $\theta$-component
from all other components.
%
\begin{lemma}[(Separation)]
\label{lemmaseparation}
There exists $C_\Delta> 0$ such that
for $\theta' \neq\theta$
\[
\E_l[\ucal_\theta| N=\theta]
> C_\Delta
\]
and
\[
\E_l[\ucal_\theta| N=\theta']
< C_\Delta.
\]
\end{lemma}
\begin{pf}
By Definition~\ref{defsparsification},
all $(a,b)\in\Upsilon''_\theta$ are
$T_{\theta'}$-stretched. Hence
\[
\E_l[\ucal_\theta| N=\theta]
\geq e^{-C_c}\vadjust{\goodbreak}
\]
and
\[
\E_l[\ucal_\theta| N=\theta']
= \frac{1}{\poly(\log n)}
\]
by Lemma~\ref{lemmaexpectations}.
Taking $C_\Delta= \frac{1}{2}e^{-C_c}$
gives the result.
\end{pf}
%
\begin{lemma}[(Concentration of $\ucal_\theta$)]
\label{lemmaconcentration}
For all $\eps_\ucal> 0$ and $C_\ucal> 0$,
there are $C_f > 0$,
$C_{st} > C_f$ and
$C_{sp}^- > 0$ large enough such that
for all $\theta, \theta'$ (possibly equal)
\[
\P_l\bigl[\bigl|\ucal_{\theta'} -
\E_l[\ucal_{\theta'} |N = \theta]\bigr| \geq\eps_\ucal|N = \theta
\bigr]
\leq\frac{1}{n^{C_\ucal}}.
\]
\end{lemma}
\begin{pf}
Let $\ical$ be as in Lemma~\ref{lemmaindependence},
and let $\ucal^{\ical}_\theta$ be the same as
$\ucal_\theta$ with the sum
restricted to $\ical$.
From Lemmas~\ref{lemmasparsification}
and~\ref{lemmaindependence},
conditioned on $\ical$,
$\ucal_\theta^\ical$ is a normalized
sum of $\Theta(\log n)$ independent
bounded variables. Concentration
of $\ucal_\theta^\ical$ therefore
follows from Lemma~\ref{lemmaazuma}
using
$m = \Omega(\log n)$, $t = O(1/\log n)$
and $\zeta= \frac{1}{2}\eps_\ucal$.
Taking $\eps_\ical= \frac{1}{2} \eps_\ucal\max_i z_i^2$
and
$C_\ical> C_\ucal$
in Lemma~\ref{lemmaindependence}
as well as $C_{sp}^- > 0$ large enough
gives the result.
\end{pf}

\section{Constructing
the clustering statistic
from data}
\label{sectionconstructing}

In this section, we provide details on the plan
laid out in Section~\ref{sectionoverview}.

Fix a GTR matrix $Q$ and constants
$\Theta\geq2$,
$0 < f \leq g < +\infty$ and $\minnu> 0$.
Let $\lambda$ be a permutation-invariant
probability measure on $\gmassump[f,g,\minnu,n]$.
In this section, we work directly with
samples $\{\sigma_X^i\}_{i=1}^k$ generated
from an \textit{unknown}
$\Theta$-mixture model $(\mixed,\nu,Q)$ picked
according to $\l$.

Our goal is to construct the clustering statistics
$\{\ucal_\theta\}_{\theta=1}^\Theta$ from
$\{\sigma_X^i\}_{i=1}^k$. These statistics
will be used in the next section to reconstruct
the topologies of the model $(\mixed,\nu,Q)$.

\subsection{Clustering algorithm}

We proceed in three steps.
Let
\[
C_c
= -\ln\biggl( \frac{1}{3 \Theta(1 - \minnu)} \minnu e^{-4g} \biggr)
\]
and
\[
\omega
= \tfrac{2}{3}\minnu e^{-4g}.
\]
The algorithm is the following:
\begin{longlist}[(2)]
\item[(1)] (\textit{Finding quasicherries})
For all pairs of leaves $a,b \in[n]$,
compute $\hat q(a,b)$, and set
\[
\hatupsilon' = \{(a,b) \in\nsqneq\dvtx\hat q(a,b) \geq\omega\}.
\]

\item[(2)] (\textit{Sparsification})
Construct $\hatupsilon''$ by keeping
each $(a,b) \in\hatupsilon'$ independently
with probability
\[
p_{sp} = \frac{C_{sp}^p \log n}{n}.\vadjust{\goodbreak}
\]

\item[(3)] (\textit{Inferring clusters})
For all $c_1 \neq c_2 \in\hatupsilon'$,
compute $\hat r(c_1,c_2)$, and
set $c_1 \sim c_2$ if
\[
\hat r(c_1,c_2)
> 0.
\]
Let $\hatupsilon''_\theta$, $\theta=1,\ldots,
\widehat{\Theta}$,
be the equivalence classes
of the transitive closure of~$\sim$.

\item[(4)] (\textit{Final sets})
Return $\hatupsilon''_{\theta}$, $\theta\in[\Theta]$.

\end{longlist}

\subsection{Analysis of the clustering algorithm}

We show that each step of the previous algorithm
succeeds with high probability.
%
\begin{lemma}[(Finding quasicherries)]
\label{lemmaquasicherries}
The set $\hatupsilon'$
satisfies the following, except with
probability at most $\exp(-{\poly(n)})$
under $\P_l$:
\[
\Upsilon_{4g} \subseteq\hatupsilon' \subseteq\Upsilon_{C_c}.
\]
\end{lemma}
\begin{pf}
We prove both inclusions.
For all $\theta\in[\Theta]$ and
$(a,b) \in\Upsilon_{4g,\theta}$,
\[
q_\theta(a,b) \geq e^{-4g}
\]
and
\[
q(a,b) \geq\minnu e^{-4g} > \tfrac{2}{3}\minnu e^{-4g}
= \omega.
\]
By Lemma~\ref{lemmahatq},
\[
\hat q(a,b) \geq\omega,
\]
except with probability $\exp(-{\poly(n)})$.\vspace*{2pt}

Similarly for any $(a,b) \in\hatupsilon'$,
by Lemma~\ref{lemmahatq},
if
\[
\hat q(a,b) \geq\omega= \tfrac{2}{3}\minnu e^{-4g},
\]
then
\[
q(a,b) \geq\tfrac{1}{3} \minnu e^{-4g},
\]
so that there is $\theta\in[\Theta]$
with
\[
\nu_\theta q_\theta(a,b) \geq\frac{1}{3 \Theta} \minnu e^{-4g}.
\]
That is,
\[
q_\theta(a,b) \geq\frac{1}{3 \Theta(1 - \minnu)} \minnu e^{-4g}
\]
and
\[
d_\theta(a,b) \leq
-\ln\biggl( \frac{1}{3 \Theta(1 - \minnu)} \minnu e^{-4g} \biggr)
= C_c.
\]
Hence $(a,b) \in\Upsilon_{C_c,\theta}$.
\end{pf}
%
\begin{lemma}[(Sparsification)]
\label{lemmasparsification2}
Assuming that the conclusions of
Lem\-ma~\ref{lemmaquasicherries} hold,
$\hatupsilon''$ is properly sparse, except with probability
$1/\poly(n)$.\vadjust{\goodbreak}
\end{lemma}
\begin{pf}
This follows from
Lemma~\ref{lemmaquasicherries}
and the choice of $p_{sp}$.
\end{pf}
%
\begin{lemma}[(Inferring clusters)]
\label{lemmaclusters}
Assuming that the conclusions of
Lemmas~\ref{lemmaquasicherries}
and~\ref{lemmasparsification2} hold,
we have $\widehat{\Theta} = \Theta$, and
there is a bijective mapping $h$ of $[\Theta]$
such that
\[
\hatupsilon''_{h(\theta)} = \Upsilon''_{\theta}
\]
with the choice $\Upsilon' = \hatupsilon'$
in Section~\ref{seccombinatorial},
except with probability $\exp(-{\poly(n)})$.
\end{lemma}
\begin{pf}
It follows from Lemmas~\ref{lemmahatr}
and~\ref{lemmaexpectations} that
$\sim$ is an equivalence relation
with equivalence classes $\Upsilon''_\theta$,
$\theta=1,\ldots,\Theta$, except
with probability $\exp(-{\poly(n)})$.
\end{pf}

\section{Tree reconstruction}
\label{sectiontree}

We now show how to use the clustering
statistics to build the topologies.
The algorithm is composed of
two steps: we first bin the sites according
to the value of the clustering statistics;
we then use the sites in
one of those bins and apply a standard distance-based
reconstruction method.
We show that the content
of the bins is made of
sites from the same component---thus
reducing the situation to the unmixed case.

Let
\begin{eqnarray*}
C_\Delta
&=& \tfrac{1}{2}e^{-C_c},
\\
\eps_\ucal&=& \tfrac{1}{3} e^{-C_c}
\end{eqnarray*}
and
\[
\eps_\ical= \frac{1}{2} \eps_\ucal\max_i z_i^2.
\]
Moreover take $C_f$, $C_{st}$, $C_{sp}^p$
and $C_{sp}^-$
so that the lemmas in Section~\ref{sectionpreliminaries}
hold.

To simplify notation, we rename the
components so that $h$ is the \textit{identity}.

\subsection{Site binning}

Let $\hatupsilon''_{\theta}$, $\theta\in[\Theta]$,
be the sets returned by the algorithm
in Section~\ref{sectionconstructing}.
Assume that the conclusions of
Lemmas~\ref{lemmaquasicherries},
\ref{lemmasparsification2}
and~\ref{lemmaclusters} hold.
We bin the sites with the following procedure:
\begin{longlist}[(2)]
\item[(1)] (\textit{Clustering statistics})
For all $i = 1, \ldots, k$
and all $\theta= 1,\ldots,\Theta$,
compute
\[
\hatucal_\theta^i
=
\frac{1}{|\hatupsilon''_\theta|}
\sum_{(a,b) \in\hatupsilon''_\theta}
\sigma_a^i \sigma_b^i.
\]

\item[(2)] (\textit{Binning sites})
For all $\theta= 1,\ldots,\Theta$,
set
\[
\hatkcal_\theta
= \{
i \in[k]\dvtx
\hatucal_\theta^i > C_\Delta
\}.
\]
\end{longlist}
We show that the binning is successful
with high probability.
%
\begin{lemma}[(Binning the sites)]
\label{lemmabinning}
Assume that the conclusions of
Lemmas~\ref{lemmaquasicherries},
\ref{lemmasparsification2}
and~\ref{lemmaclusters} hold.
For any $C_k$, there exists $C_\ucal$
large enough so that, for all
$\theta\in[\Theta]$,
\[
\hatkcal_{\theta}
= \kcal_\theta,
\]
except with probability $1/\poly(n)$.
\end{lemma}
\begin{pf}
This follows from Lemmas~\ref{lemmaseparation}
and~\ref{lemmaconcentration} by a union
bound over all samples.
\end{pf}

\subsection{Estimating a distorted metric}
\label{secdistorted}

\subsubsection*{Estimating evolutionary distances}
We estimate evolutionary distances on each component.
For all $\theta\in[\Theta]$,
let $\hatkcal_{\theta}$ be as above and assume
the conclusions of Lemma~\ref{lemmabinning}
hold.
\begin{longlist}[(2)]
\item[(1)] (\textit{Estimating distances})
For all $\theta= 1,\ldots,\Theta$
and $a\neq b \in[n]$,
compute
\[
\hat q_\theta(a,b) =
\frac{1}{|\hatkcal_\theta|}
\sum_{i\in\hatkcal_\theta} \sigma^i_a \sigma^i_b.
\]
\end{longlist}

\begin{lemma}[(Estimating distances)]
\label{lemmadistances}
Assume the conclusions of
Lem\-ma~\ref{lemmabinning} hold.
The following hold
except with probability $\exp(-{\poly(n)})$:
for all $\theta\in[\Theta]$
and all $a\neq b \in[n]$,
\[
|
\hat q_\theta(a,b) - q_\theta(a,b)
|
\leq
\frac{1}{n^{C_q}}.
\]
\end{lemma}
\begin{pf}
The result follows from
Lemma~\ref{lemmahatq}.
\end{pf}

\subsubsection*{Tree construction}
To reconstruct the tree,
we use a distance-based method
of~\cite{DaMoRo09}. We require the following
definition.
%
\begin{definition}[(Distorted metric~\cite{Mossel07,KiZhZh03})]
\label{defdistortedmetric}
Let $T = (V,E;\phi;w)$ be a phylogeny
with corresponding tree metric $d$, and let $\tau, \Psi> 0$.
We say that $\hat d\dvtx X\times X \to(0,+\infty]$ is a $(\tau, \Psi
)$-\textit{distorted metric}
for $T$ or a $(\tau, \Psi)$-\textit{distortion} of $\dist$ if:
\begin{longlist}[(2)]
\item[(1)] (\textit{Symmetry})
For all $a, b \in X$, $\hat d$ is symmetric, that is,
\[
\hat d(a,b) = \hat d(b,a);
\]
\item[(2)] (\textit{Distortion}) $\hat d$ is accurate on
``short'' distances; that is, for all $a, b \in X$, if either
$d (a,b) < \Psi+ \tau$ or $\hat d(a,b) < \Psi+ \tau$, then
\[
|d(a,b) - \hat d(a,b)| < \tau.
\]
\end{longlist}
\end{definition}

An immediate consequence of~\cite{DaMoRo09}, Theorem 1,
is the following.\vadjust{\goodbreak}
%
\begin{uclaim*}[(Reconstruction from distorted metrics~\cite{DaMoRo09})]
Let $T = (V,E;\phi;w)$ be a phylogeny
in
$\phy_{f,g}$.
Then the topology of $T$ can be recovered
in polynomial time from a
$(\tau,\Psi)$-distortion $\hat d$ of $d$ as long as
\[
\tau\leq\frac{f}{5}
\]
and
\[
\Psi\geq5 g \log n.
\]
\end{uclaim*}
%
\begin{remark}
The constants above are not optimal
but will suffice for our purposes.
\end{remark}

See~\cite{DaMoRo09} for the details
of the reconstruction algorithm.

We now show how to obtain a
$(f/5,5g\log n)$-distortion with high
probability for each component.
%
\begin{lemma}[(Distortion estimation)]
\label{lemmadistortion}
There exist $C_q,C_k > 0$
so that, given that the conclusions
of Lemma~\ref{lemmadistances} hold,
for all $\theta\in[\Theta]$,
\[
\hat d_\theta(a,b) = -\ln( \hat{q}_\theta(a,b)_+ ),\qquad
(a,b) \in X\times X,
\]
is a $(f/5,5 g\log n)$-distortion of
$d_\theta$.
\end{lemma}
\begin{pf}
Fix $\theta\in[\Theta]$.
Define
\[
\short= \{(a,b) \in X \times X\dvtx
d_\theta(a,b) \leq15 g \log n\}
\]
and
\[
\notshort= \{(a,b) \in X \times X\dvtx
d_\theta(a,b) > 12 g \log n\}.
\]

Let $(a,b)\in\short$. Note that
\[
e^{-d_\theta(a,b)}
\geq\exp(-
15 g \log n)
\equiv\frac{1}{n^{C_q'}},
\]
where the last equality is a definition.
Then, taking $C_q$ (and hence $C_k$) large enough,
from Lemma~\ref{lemmadistances}, we have
\[
|
\hat d_\theta(a,b)
- d_\theta(a,b)
|
\leq
\frac{f}{5}.
\]

Similarly,
let $(a,b)\in\notshort$. Note that
\[
e^{- d_\theta(a,b)}
< \exp(-12 g \log n)
\equiv\frac{1}{n^{C_q''}},
\]
where the last equality is a definition. Then, taking $C_q$ large
enough, from Lem\-ma~\ref{lemmadistances} we have
\[
\hat d_\theta(a,b)
\geq
5 g \log n + \frac{f}{5}.
\]
\upqed\end{pf}

\section{Proof of main theorems}

We are now ready to prove the main theorems.
\begin{pf*}{Proof of Theorem~\ref{thm3}}
Let $C_1, C_2 > 0$.
Let $A_n$ be the subset of
those $\Theta$-mixture models $(\mixed,\nu,Q)$ in
$\gmassump[f,g,\minnu,n]$
for which
part (3) of Lemma~\ref{lemmasparsification}
holds with probability at least $1 - n^{-C_1}$
under the random choices of the algorithm.
By the proof of Lemma~\ref{lemmasparsification},
for small enough $C_1, C_2 > 0$,
we have $\lambda_n[A_n^c] \leq n^{-C_2}$.
On $A_n$, the lemmas in
Sections~\ref{sectionpreliminaries},
\ref{sectionconstructing}
and
\ref{sectiontree}
hold with probability $1 - 1/\poly(n)$.
Then the topologies are correctly reconstructed
by the claim in Section~\ref{secdistorted}.
\end{pf*}
\begin{pf*}{Proof of Theorem~\ref{thm1}}
Let
\[
(\mixed,\nu,Q) \nsim(\mixed',\nu',Q)
\in\bigcup_{n \geq1} A_n.
\]
Then, by Theorem~\ref{thm3},
the algorithm correctly reconstructs
the topologies in $(\mixed,\nu,Q)$ with
probability $1 - 1/\poly(n)$ on
sequences of length $k = \poly(n)$.
Repeating the reconstruction on
independent sequences and
taking a majority vote, we get
almost sure convergence to
the correct topologies. The same
holds for $(\mixed',\nu',Q)$.
Hence,
\[
\dcal_l[(\mixed,\nu,Q)]
\neq\dcal_l[(\mixed',\nu',Q)].
\]
\upqed\end{pf*}
\begin{pf*}{Proof of Theorem~\ref{thm2}}
Let
\[
(\mixed,\nu,Q) \in\bigcup_{n \geq1} A_n
\]
with $\Theta= 2$ and $\nu= (1/2,1/2)$.
Then,
from the proof of Lemma~\ref{lemmabinning},
there exists a clustering statistic
such that samples from $T_1$ and $T_2$
are correctly distinguished
with probability $1 - 1/\poly(n)$.
Recall that
\[
\tv{\dcal- \dcal'}
= \sup_{B \in\fcal} |
\dcal(B) - \dcal'(B)
|.
\]
Taking $B$ to be the event that
a site is recognized as belonging
to component $1$ by the
clustering statistic above, we get
\[
\tv{\dcal_l[T_1,Q] - \dcal_l[T_2,Q]}
= 1 - o_n(1).
\]
\upqed\end{pf*}

\section{Concluding remarks}

Our techniques also admit the following
extensions:
\begin{itemize}
\item When $Q$ is unknown, one
can still apply our technique by using
the following idea.
Note that all we need is an eigenvector
of $Q$ with negative eigenvalue.
Choose a pair $(a,b)$ of close
leaves using, for instance,
the classical log-det distance~\cite{SempleSteel03}.
Under a permutation-invariant
measure, $(a,b)$ is stretched in all
but one component, with high probability.
One can then compute an eigenvector decomposition
of the transition matrix between $a$ and $b$.
We leave out the details.\vadjust{\goodbreak}

\item The minimum
frequency assumption is not necessary
as long as one has an upper bound
on the number of components
and that one requires only that
frequent enough components
be detected and reconstructed.
We leave out the details.
\end{itemize}



\printaddresses

\end{document}